\documentclass{amsart}

\usepackage{amsfonts,amsmath,amssymb,latexsym}

\numberwithin{equation}{section}
\newtheorem{theorem}{Theorem}[section]

\newtheorem{corollary}[theorem]{Corollary}
\newtheorem{definition}[theorem]{Definition}

\newtheorem{example}[theorem]{Example}
\newtheorem{notation}[theorem]{Notation}

\begin{document}

\pagenumbering{arabic} \pagestyle{headings}
\def\sof{\hfill\rule{2mm}{2mm}}
\def\SS{\mathcal S}
\def\qq{{\bold q}}
\def\txx{{\frac1{2\sqrt{x}}}}
\def\mn{\text{-}}
\def\gp{\text{Gp}}

\title{{\sc Counting occurrences of some subword patterns}}

\author{Alexander Burstein \and Toufik Mansour}

\date{January 26, 2003}

\address{Department of Mathematics, Iowa State University, Ames, IA
50011-2064 USA}

\email{burstein@math.iastate.edu}

\address{LaBRI, Universit\'e Bordeaux 1, 351 cours de la
Lib\'eration, 33405 Talence Cedex, France}

\email{toufik@labri.fr}

\begin{abstract}
We find generating functions for the number of strings (words)
containing a specified number of occurrences of certain types of
order-isomorphic classes of substrings called subword patterns. In
particular, we find generating functions for the number of strings
containing a specified number of occurrences of a given $3$-letter
subword pattern.
\end{abstract}

\maketitle

\section{Introduction}

Counting the number of words which contain a set of given strings
as substrings a certain number of times is a classical problem in
combinatorics. This problem can, for example, 
be attacked using the transfer matrix
method (see \cite[Section 4.7]{St}). In particular, 
it is a well-known fact that the generating function of such words is always
rational. For example, in \cite[Example 4.7.5]{St} it is shown
that the generating function for the number of words in $[3]^n$ where
neither $11$ nor $23$ appear as two consecutive digits is given by
\[
\frac{3+x-x^2}{1-2x-x^2+x^3}.
\]

In this paper, we present,  in several cases, 
a complete solution for the problem of the enumeration
of words containing a {\em subword pattern} (see below for the
precise definition) of length $l$ exactly $r$
times. 
For example, we find the number of words in $[3]^n$
containing the subword pattern $111$ exactly $r$ times, that is,
the number of words which contain $111$, $222$, and $333$ as
substrings a total of $r$ times. 

R\'egnier and Szpankowski \cite{RS} used a combinatorial approach
to study the frequency of occurrences of strings 
(which they also called a ``pattern'') from a given
set in a random word, when
overlapping copies of the ``patterns'' are counted separately (see
\cite[Theorem 2.1]{RS}). We note that the term ``pattern'' in
\cite{RS} is used to denote an exact string rather than its type
with respect to order isomorphism. For example, the ``pattern'' $112$
in \cite{RS} is the actual string $112$, whereas in our setting an
occurrence of the
subword pattern $112$ is any substring $aab$ of the ambient string
with $a<b$. Although, in principle, it is possible to deduce our
results from the result by R\'egnier and Szpankowski, our direct
derivations are much simpler.

Goulden and Jackson~\cite{GJ} also consider sequences with
distinguished substrings and use the term ``pattern of a
sequence''. However, their ``pattern'' is more locally defined
than in this paper in that only order relations between adjacent
elements of a string are considered, rather than order relations
between any pair of elements of a string, as is done in this
paper. For example, the pattern ``rise, non-rise'' (or $(<,\ge)$, or
$\pi_1\pi_2$) as defined in \cite{GJ} includes the subword patterns
$121, 122, 132, 231$ as defined in this paper. However, we show
that each of the subword patterns $121, 122, 132$ is avoided by a
different number of words (of a given length on a given alphabet)
than the other two patterns.

In what follows, we use analytical and combinatorial means to find
a complete answer for several cases of counting strings with a
specified number of occurrences of certain patterns.

\subsection{Classical patterns in permutations} Let $\pi\in S_n$
and $\tau\in S_m$ be two permutations. An \emph{occurrence\/} of
$\tau$ in $\pi$ is a subsequence $1\le i_1<i_2<\dots<i_m\le n$
such that $(\pi(i_1),\dots,\pi(i_m))$ is order-isomorphic to
$\tau$. In this context, $\tau$ is usually called a
\emph{pattern}. We denote the number of occurrences of $\tau$ in
$\pi$ by $\pi(\tau)$.

Recently, much attention has been paid to the problem of counting
the number of permutations of length $n$ containing a given number
$r\ge0$ of occurrences of a certain pattern $\tau$. Most of the
authors consider only the case $r=0$, thus studying permutations
\emph{avoiding\/} a given pattern. Only a few papers consider the
case $r>0$, usually restricting themselves to the patterns of
length $3$. In fact, simple algebraic considerations show that
there are only two essentially different cases for $\tau\in S_3$,
namely, $\tau=123$ and $\tau=132$. Noonan \cite{No} has proved
that the number of permutations in $S_n$ containing $123$ exactly
once is given by $\frac 3n\binom{2n}{n-3}$. A general approach to
the problem was suggested by Noonan and Zeilberger \cite{NZ}; they
gave another proof of Noonan's result, and conjectured that the
number of permutations in $S_n$ containing $123$ exactly twice is
given by $\frac{59n^2+117n+100}{2n(2n-1)(n+5)}\binom{2n}{n-4}$ and
the number of permutations in $S_n$ containing $132$ exactly once
is given by $\binom{2n-3}{n-3}$. The first conjecture was proved
by Fulmek \cite{Fu} and the second conjecture was proved by B\'ona
in \cite{B2}. A general conjecture of Noonan and Zeilberger states
that the number of permutations in $S_n$ containing $\tau$ exactly
$r$ times is $P$-recursive in $n$ for any $r$ and $\tau$. It was
proved by B\'ona \cite{B1} for $\tau=132$.  However, as stated in
\cite{B1}, a challenging question is to describe the number of
permutations in $S_n$ containing $\tau\in S_3$ exactly $r$ times,
explicitly for any given $r$. Later, Mansour and Vainshtein
\cite{MVr} suggested a new approach to this problem in the case
$\tau=132$, which allows one to get an explicit expression for the
number of permutations in $S_n$ containing $132$ exactly $r$ times
for any given $r$.


\subsection{Generalized patterns in permutations} In \cite{BS},
Babson and Steingr\'{\i}msson introduced generalized permutation
patterns that allow the requirement that two adjacent letters in a
pattern must be adjacent in the permutation. For example, an
occurrence of a generalized pattern $12\mn3$ in a permutation $\pi
= a_1 a_2 \cdots a_n$ is a subword $a_i a_{i+1} a_{j}$ of $\pi$
such that $a_i<a_{i+1}<a_{j}$.

\begin{notation}
\rm Unfortunately, there is a bit of confusion in denoting
classical and generalized patterns. Before generalized patterns
were introduced, the hyphens were unnecessary, hence classical
patterns (those with all possible hyphens) are often written with
no hyphens when generalized patterns are not considered, and with
all hyphens when they are. For example, a classical pattern $123$
is now denoted by $1\mn2\mn3$ when considered as a generalized
pattern (using the notation of \cite{BS, Claesson}). Unless
otherwise stated, all patterns under consideration from now on are
generalized patterns.
\end{notation}

In \cite{EN} Elizalde and Noy presented the following theorem
regarding the distribution of the number of occurrences of any
generalized pattern of length $3$ without hyphens.

\begin{theorem}{\em(Elizalde and Noy \cite{EN})}
  Let $h(x)=\sqrt{(x-1)(x+3)}$. Then
  \begin{align*}
    \sum_{n\ge 0}\sum_{\pi\in S_n} x^{\pi(123)}\frac{t^n}{n!}
    &=\frac{2h(x)e^{\frac{1}{2}(h(x)-x+1)t}}
    {h(x)+x+1+(h(x)-x-1)e^{h(x)t}},\\
    \sum_{n\ge0}\sum_{\pi\in S_n} x^{\pi(213)}\frac{t^n}{n!}
    &=\frac{1}{1-\int_0^{t}e^{(x-1)z^2/2}dz},
  \end{align*}
where $\pi(123)$ (respectively, $\pi(213)$) is the number of
occurrences of the \emph{generalized} pattern $123$ (respectively,
$213$) without hyphens in $\pi$.
\end{theorem}

On the other hand, Claesson \cite{Claesson} gave a complete answer
for the number of permutations avoiding a generalized pattern of
the form $xy\mn z$ where $xyz\in S_3$. Later, Claesson and Mansour
\cite{CM} presented an algorithm to count the number of
permutations containing a generalized pattern of the form $xy\mn
z$ exactly $r$ times for any given $r\ge0$, where $xyz\in S_3$.

\begin{theorem}{\rm(Claesson and Mansour \cite{CM})}
The ordinary generating function for the number of permutations of
length $n$ avoiding the generalized pattern $12\mn3$ (or $23\mn1$)
is
\[
\sum_{k\ge 0} \frac{x^k}{(1-x)(1-2x)\cdots(1-kx)}.
\]

The ordinary generating function for the number of permutations of
length $n$ avoiding  the generalized pattern $2\mn13$ is
\[
C(x)=\frac{1-\sqrt{1-4x}}{2x}.
\]

The ordinary generating function for the number of permutations of
length $n$ containing exactly one occurrence of the generalized
pattern $12\mn3$ is
\[
\sum_{n\ge 1}\frac x {1-nx}\sum_{k\ge 0} \frac {kx^{k+n}}
{(1-x)(1-2x)\cdots(1-(k+n)x)}.
\]

The ordinary generating function for the number of permutations of
length $n$ containing exactly one occurrence of the generalized
pattern  $23\mn1$ is
\[
\sum_{n\ge 1}\frac x {1-(n-1)x}\sum_{k\ge 0} \frac {kx^{k+n}}
{(1-x)(1-2x)\cdots(1-(k+n)x)}.
\]

The ordinary generating function for the number of permutations of
length $n$ containing exactly one occurrence of the generalized
pattern $2\mn13$ is
\[
\frac{x^3C(x)^7}{1-tC(x)^2}.
\]
\end{theorem}
\subsection{Generalized patterns in words} A \emph{generalized
pattern} $\tau$ is a (possibly hyphenated) string in $[\ell]^m$
which contains all letters from $[\ell]=\{1,\dots,\ell\}$. We say
that the string $\sigma\in[k]^n$ \emph{contains} a generalized
pattern $\tau$ exactly $r$ times (denoted by $r=\sigma(\tau)$) if
$\sigma$ contains $r$ different subsequences isomorphic to $\tau$
in which the entries corresponding to consecutive entries of
$\tau$ not separated by a hyphen must be adjacent. We call the
generalized patterns without hyphens \emph{subword patterns}. If
$r=0$, we say that $\sigma$ \emph{avoids} $\tau$ and write
$\sigma\in [k]^n(\tau)$. Thus, $[k]^n(\tau)$ denotes the set of
strings in $[k]^n$ (i.e., $n$-long $k$-ary strings) which avoid
$\tau$. For example, a string $\pi=a_1a_2\dots a_n$ avoids the
generalized pattern $12\mn1$ if $\pi$ has no subsequence
$a_ia_{i+1}a_j$ with $j>i+1$ and $a_i=a_j<a_{i+1}$.

\begin{example}
\rm Davenport-Schinzel sequences \cite{DS} can be defined in terms
of subword pattern avoidance as follows. For any $d\geq 1$, let
$T_d$ be the set of all the subword patterns $\pi=a_1a_2\cdots
a_{d+1}\in [d+1]^{d+1}$ such that either
$a_{2j}<a_{2j+1}>a_{2j+2}$ for all $j$, or
$a_{2j-1}<a_{2j}>a_{2j+1}$ for all $j$. For example, $T_2=\{121,
132, 231, 212, 213, 312\}$. An $k$-ary $n$-long sequence avoiding
the subword pattern $11$ (i.e., with no equal consecutive letters)
and avoiding all the subword patterns in $T_d$ (there are no
alternating subwords of length greater than $d+1$) is called a
Davenport-Schinzel sequence if $n$ is maximal.
\end{example}

Let $f_{\tau;r}(n,k)$ be the number of words $\sigma\in[k]^n$ such
that $\sigma(\tau)=r$. Denote the corresponding bivariate
generating function by $F_{\tau}(x,y;k)$, in other words,
\[
F_{\tau}(x,y;k)=\sum_{n\ge0}\sum_{r\ge0}{f_{\tau;r}(n,k)x^ny^r}.
\]

Burstein \cite{Burstein} gave a complete answer for the numbers
$f_{\tau;0}(n,k)$ where $\tau$ is a $3$-letter classical pattern.
Later, Burstein and Mansour \cite{BM1,BM2} presented a complete
answer for the number $f_{\tau;0}(n,k)$ where $\tau$ is a generalized
pattern of length $3$ (a word of length $3$).

In this paper, we present a complete answer for several cases of
$f_{\tau;r}(n,k)$ where $\tau$ is a subword pattern of length $l$
(which is the analogue of the results by Elizalde and Noy in \cite{EN}).
In particular, we find a complete answer for the case $l=3$.

\section{Counting a subword pattern of length $l$}
In this section we find $F_\tau(x,y;k)$ for several cases of
$\tau$. Burstein and Mansour \cite{BM2} found $F_\tau(x,y;k)$ for
the subword pattern $\tau=11\dots1\in[1]^l$ and proved the 
following theorem.

\begin{theorem}\label{tones}{\rm(Burstein and Mansour
\cite[Th. 2.1]{BM2})} Let $\tau=11\dots 1\in [1]^l$ be a subword
pattern. Then
\[
F_{\tau}(x,y;k)=\frac{1+(1-y)x\sum_{j=0}^{l-2} (kx)^j
-(1-y)(k-1)\sum_{d=2}^{l-1}x^d\sum_{j=0}^{l-1-d}
(kx)^j}{1-(k-1+y)x-(k-1)(1-y)(1-x^{l-2})\frac{x^2}{1-x}}.
\]
\end{theorem}

\subsection{The subword pattern $\tau=11\dots12$} Let $\tau=11\dots12\in[2]^l$ be a
subword pattern. Define $d_{\tau;r}(n,k)$ to be the number of
words $\beta\in[k]^n$ such that $(\beta,k+1)$ contains $\tau$
exactly $r$ times, and denote the corresponding generating
function by $D_{\tau}(x,y;k)=\sum\limits_{n,r\geq0}
d_{\tau;r}(n,k)x^ny^r$. Let us find a recurrence for $F_{\tau}$.

Let $\sigma\in [k]^n$ ($k\ge 2$), with $\sigma(\tau)=r$, contain
exactly $d$ copies of the letter $k$. If $d=0$, then
$\sigma\in[k-1]^n$ and $\sigma(\tau)=r$. If $d\ge 1$, then
$\sigma=\sigma_1 k\sigma_2$, where $\sigma_1\in [k-1]^{n_1}$,
$\sigma_2\in [k]^{n_2}$, $n_1+n_2+1=n$ and
$(\sigma_1,k)(\tau)+\sigma_2(\tau)=r$. Taking generating
functions, we see that the above translates into
\[
F_\tau(x,y;k)=F_\tau(x,y;k-1)+xD_\tau(x,y;k-1)F_\tau(x,y;k),
\]
or, equivalently,
\[
F_{\tau}(x,y;k)=\frac{F_{\tau}(x,y;k-1)}{1-xD_{\tau}(x,y;k-1)}.
\]

Let us now find the recurrence for $D_\tau$. Let $\sigma\in[k]^n$
be such that $(\sigma,k+1)$ contains $\tau$ exactly $r$ times, and
has exactly $d$ letters $k$. Then $\sigma=\sigma_1 k\sigma_2
k\dots k\sigma_d k\sigma_{d+1}$ for some $\sigma_i\in
[k-1]^{n_i}$, $1\le i\le d+1$, where $n_1+\dots+n_{d+1}=n-d$ and
$$(\sigma_1,k)(\tau)+\dots+(\sigma_d,k)(\tau)+(\sigma_{d+1},k+1)(\tau)+\delta(\sigma
\text{ ends on } l-1\ k\text{'s})=r.$$ Taking generating functions,
we obtain
\begin{multline*}
D_\tau(x,y;k)=\sum_{d=0}^{l-2}{x^dD_\tau^{d+1}(x,y;k-1)}\\
+\sum_{d=l-1}^{\infty}x^d\Big(D_\tau^{d+1}(x,y;k-1)-D_\tau^{(d+1)-(l-1)}(x,y;k-1)\\
+yD_\tau^{(d+1)-(l-1)}(x,y;k-1)\Big),
\end{multline*}
which, after summing over $d$, yields
\[
D_\tau(x,y;k-1)=\frac{(1-x^{l-1}(1-y))D_\tau(x,y;k-1)}{1-xD_\tau(x,y;k-1)}.
\]
These two recurrences, together with
$D_{\tau}(x,y;0)=F_\tau(x,y;0)=1$ and induction on $k$, yield the
following theorem.

\begin{theorem}\label{th112}
Let $\tau=11\dots12\in[2]^l$ be a subword pattern such that $l\ge
1$; then
\[
F_\tau(x,y;k)=\frac{1-y}{1-x^{2-l}-y+x^{2-l}(1-x^{l-1}(1-y))^k}.
\]
\end{theorem}

\begin{example}\label{ex112}
\rm (see Burstein and Mansour \cite[Th. 3.10]{BM2}) Letting $l=3$
and $y=0$ in Theorem \ref{th112}, we get that the generating
function for the number of words in $[k]^n$ avoiding the subword
pattern $112$ is given by
\[
\frac{1}{1-\frac{1}{x}+\frac{1}{x}(1-x^2)^k}.
\]
\end{example}

In the special case of $l=3$, we get from Theorem \ref{th112} the
following result.

\begin{corollary}
The generating function for the number of words in $[2]^n$
containing the subword pattern $112$ exactly $r$ times is given by
\[
\frac{x^{3r}}{(1-x)^{r+1}(1-x-x^2)^{r+1}}.
\]
\end{corollary}
\begin{proof}
Let $\tau=112$ be a subword pattern. It is easy to see that a word
$\sigma\in[2]^n$ with $\sigma(\tau)=r$ must have the form
$\sigma=\sigma_1\tau\sigma_2\tau\dots\tau\sigma_{r+1}$, for some
$\sigma_1,\dots,\sigma_{r+1}\in[2]^n(\tau)$. Now, from Example
\ref{ex112} with $k=2$, we have that
$F_{\tau}(x;2)=\frac{1}{(1-x)(1-x-x^2)}$, hence the result
follows.
\end{proof}

\subsection{The subword pattern $\tau=211\dots112$}
Let $\tau=211\dots112\in [2]^l$ be a subword pattern. We define
$d_\tau(n,r;k)$ to be the number of words $\beta\in[k]^n$ such that
$(k+1,\beta,k+1)$ contains $\tau$ exactly $r$ times, and denote
the corresponding generating function by $D_\tau(x,y;k)$. Let
$\sigma\in[k]^n$ such that $\sigma(\tau)=r$, and such that $\sigma$ 
contains $d$
occurrences of the letter $k$. For $d=0$, the generating function
for the number of such words $\sigma$ is given by
$F_\tau(x,y;k-1)$, and for $d\ge 1$, by
$x^dF_\tau^2(x,y;k-1)D_\tau^{d-1}(x,y;k-1)$ (since in that case
$\sigma=\sigma_0k\sigma_1k\dots k\sigma_{d-1}k\sigma_d$, where all
$\sigma_i\in [k-1]^{n_i}$, $\sum{n_i}=n-d$, and
$\sigma(\tau)=\sigma_0(\tau)+(k,\sigma_1,k)(\tau)+\dots+(k,\sigma_{d-1},k)(\tau)+\sigma_d(\tau)$).
Hence, if we sum over all $d\ge 0$, we get
\[
F_\tau(x,y;k)=F_\tau(x,y;k-1)+\frac{xF_\tau^2(x,y;k-1)}{1-xD_\tau(x,y;k-1)}.
\]
On the other hand, the word $(k+1,\beta,k+1)$, with $\beta$ as above, 
contains an occurrence
of $\tau$ involving the two letters $k+1$ if and only if $\beta$
is a constant string of length $l-2$, otherwise,
$(k+1,\beta,k+1)(\tau)=\beta(\tau)$. Taking generating functions,
we obtain
\[
D_\tau(x,y;k)=kx^{l-2}y+F_\tau(x,y;k)-kx^{l-2}.
\]
Therefore, using the initial conditions
$F_\tau(x,y;0)=D_\tau(x,y;0)=1$ and induction on $k$, we get the
following theorem.

\begin{theorem}\label{th212}
Let $\tau=211\dots112\in[2]^l$ be a subword pattern and $l\ge 2$,
then
\[
F_\tau(x,y;k)=\frac{1}{1-x-x\sum_{j=0}^{k}\frac{1}{1+jx^{l-1}(1-y)}}.
\]
\end{theorem}

\begin{example}\rm (Burstein and Mansour \cite[Th. 3.12]{BM2})
Letting $l=3$ and $y=0$ in Theorem \ref{th212}, we get that the
generating function for the number of words in $[k]^n$ avoiding
the subword pattern $212$ is given by
\[
\frac{1}{1-x-x\sum_{j=0}^{k-1}\frac{1}{1+jx^2}}.
\]
\end{example}

\subsection{The subword pattern $\tau=m\tau'm$}

\begin{theorem}\label{thll}
Let $\tau=m\tau'm\in [m]^l$ be a subword pattern, where $\tau'$
does not contain $m$. Then for $k\ge m$,
\[
F_\tau(x,y;k)=\frac{1}{1-(m-1)x-x\sum_{j=m-1}^{k-1}
\frac{1}{1+\binom{j}{m-1}x^{l-1}(1-y)}}.
\]
\end{theorem}

\begin{proof}
Let $\sigma\in [k]^n$. The generating function for the number of
words $\sigma$ which do not contain $m$ and contain $\tau$ exactly
$r$ times is given by $F_\tau(x,y;k-1)$. Now assume that the leftmost
$m$ in $\sigma$ is at position $i$. Then
$\sigma=\sigma_1m\sigma_2$ and
$\sigma(\tau)=\sigma_1(\tau)+(m,\sigma_2)(\tau)$, so the
generating function for the number of such words $\sigma$ is given
by $xF_\tau(x,y;k-1)D_\tau(x,y;k)$, where $D_\tau(x,y;k)$ is the
generating function for the number of words $\sigma\in [k]^n$ such
that $(m,\sigma)$ contains $\tau$ exactly $r$ times. Therefore,
\[
F_\tau(x,y;k)=F_\tau(x,y;k-1)+xF_\tau(x,y;k-1)D_\tau(x,y;k).
\]
On the other hand, let $\sigma'=(m,\sigma)\in [k]^{n+1}$. If
$\sigma$ does not contain $m$, then the generating function for
the number of such $\sigma$ is given by $F_\tau(x,y;k-1)$.
Otherwise, let $i$ be the position of the leftmost letter $m$ and
let $\sigma|_i$ be the left prefix of $\sigma$ of length $i$, then
the generating function for these words is given by
$x\left(F_\tau(x,y;k-1)-x^{l-2}\binom{k-1}{m-1}\right)D_\tau(x,y;k)$
if $(m,\sigma|_i)$ is not order-isomorphic to $\tau$, or by
$xyx^{l-2}\binom{k-1}{m-1}D_\tau(x,y;k)$ if $(m,\sigma|_i)$ is
order-isomorphic to $\tau$. Therefore,
\[
\begin{split}
D_\tau(x,y;k)&=F_\tau(x,y;k-1)+x\left(F_\tau(x,y;k-1)-x^{l-2}\binom{k-1}{m-1
}\right)D_\tau(x,y;k)\\
&\kern4cm+yx^{l-1}\binom{k-1}{m-1}D_\tau(x,y;k).
\end{split}
\]
Hence, from the above two equations, we obtain
\[
F_\tau(x,y;k)=\frac{\left(1+x^{l-1}\binom{k-1}{m-1}(1-y)\right)F_\tau(x,y;
k-1)} {1+x^{l-1}\binom{k-1}{m-1}(1-y)-xF_\tau(x,y;k-1)},
\]
so, by induction on $k$ with the initial condition
$F_\tau(x,y;m-1)=\frac{1}{1-(m-1)x}$, we get the desired result.
\end{proof}

\begin{example} \rm
Applying Theorem \ref{thll} to the subword patterns $2112$ and $3123$,
we get
\begin{align*}
F_{2112}(x,y;k)&=\frac{1}{1-x-x\sum_{j=0}^{k-1}\frac{1}{1+jx^3(1-y)}},\\
F_{3123}(x,y;k)&=\frac{1}{1-2x-x\sum_{j=2}^{k-1}\frac{1}{1+j(j-1)x^3(1-y)/2}}.
\end{align*}
\end{example}

\begin{definition}
We say that the patterns $\beta$ and $\gamma$ are \emph{strongly
Wilf-equivalent}, or are in the same \emph{strong Wilf class}, if
the number of words in $[k]^n$ containing $\beta$ exactly $r$
times is the same as the number of words in $[k]^n$ containing
$\gamma$ exactly $r$ times, for any $r\ge0$.
\end{definition}

By Theorem \ref{thll} and the symmetry operations ``reversal" and
``complement," we immediately get the following corollary.

\begin{corollary}\label{c41}
The subword patterns $1121$ and $1221$ are in the same strong Wilf
class.
\end{corollary}

\subsection{The subword pattern $\tau=m\tau'(m+1)$}
Let $\tau=m\tau'(m+1)$ be a subword pattern, where $\tau'$ does
not contain $m$ or $m+1$. Note that $\tau$ is in the same symmetry
class as $r(c(\tau))=1r(c(\tau'))2$. This case is treated in a
similarl manner as the case of $\tau=m\tau'm$. As a result, we obtain
the theorem below.

\begin{theorem}\label{thab}
Let $\tau=m\tau'(m+1)\in [m+1]^l$ be a subword pattern, where
$\tau'$ does not contain $m$ or $m+1$. Then for $k\ge m$,
\[
F_\tau(x,y;k)=\frac{1}{1-(m-1)x-x\sum_{i=m-2}^{k-2}
\prod_{j=m-2}^i \left(1-\binom{j}{m-1}x^{l-1}(1-y)\right)}.
\]
\end{theorem}

\section{Subword patterns of length $3$}
The symmetry class representatives of $3$-letter subword patterns
are $111$, $112$, $212$, $123$, $213$. In the current subsection,
we find explicit formulas for $F_\tau(x,y;k)$ for each of these
representatives $\tau$. Theorem \ref{tones} yields the answer for
the first class.

\begin{theorem}\label{the111}
Let $\tau=111$ be a subword pattern. Then, for all $k\ge 0$, we have
\[
F_{\tau}(x,y;k)=\frac{1+x(1+x)(1-y)}{1-(k-1+y)x-(k-1)(1-y)x^2}.
\]
\end{theorem}

Theorems \ref{th112} and \ref{th212} contain already the answers for the
second and the third classes, respectively. Let us summarize the
corresponding results in the theorem below.

\begin{theorem}\label{the112212}
Let $112$ and $212$ be subword patterns. For $k\ge 0$,
\[
\begin{split}
F_{112}(x,y;k)&=\frac{1-y}{1-\frac{1}{x}-y+\frac{1}{x}(1-x^2(1-y))^k},\\
F_{212}(x,y;k)&=\frac{1}{1-x-x\sum_{j=0}^{k}\frac{1}{1-jx^2(1-y)}}.
\end{split}
\]
\end{theorem}

Now let us find the generating function for the fourth class,
$F_{\tau}(x,y;k)$ where $\tau=123$ is a subword pattern. Let
$D_{\tau}(x,y;k)$ be the generating function for the number of
words $\sigma\in [k]^n$ such that $(\sigma,k+1)$ contains the
subword pattern $123$ exactly $r$ times. Suppose a word
$\sigma\in[k]^n(\tau)$ has exactly $d$ letters $k$. Then
$\sigma=\sigma_0k\sigma_1k\dots k\sigma_d$, where all
$\sigma_i\in[k-1]^n$, and any occurrence of $\tau$ in $\sigma$ must
be either in $(\sigma_i,k)$ for some $i=0,1,\dots,d-1$, or in
$\sigma_d$. Therefore, the generating function for the number of
such words $\sigma$ is $(xD_{\tau}(x,y;k-1))^dF_{\tau}(x,y;k-1)$,
so
\[
F_{\tau}(x,y;k)=\sum_{d\ge0}{(xD_{\tau}(x,y;k-1))^d
F_{\tau}(x,y;k-1)}.
\]

On the other hand, suppose $\sigma\in[k]^n$ is counted by
$D_{\tau}(x,y;k)$. Then $(\sigma,k+1)=\sigma_0k\sigma_1k\dots
k\sigma_dk+1$ (where $\sigma_i\in[k-1]^n$ for all $i$) contains the pattern
$\tau$ exactly $r$ times. If $d=0$, then $\sigma\in[k-1]^n$. If
$d\ge 1$, there are several possibilities. If
$\sigma_d\ne\emptyset$ or $\sigma_d=\sigma_{d-1}=\emptyset$, then
all occurrences of the pattern $\tau$ in $(\sigma,k+1)$ are in
$(\sigma_i,k)$ for some $i=0,1,\dots,d-1$, or in $(\sigma_d,k+1)$.
If $\sigma_d=\emptyset$ and $\sigma_{d-1}\ne\emptyset$, then there
is one extra occurrence of $\tau$ since $(\sigma,k+1)$ ends by
$(a,k,k+1)$ for some $a<k$. Taking generating functions, we obtain
\[
\begin{split}
D_{\tau}(x,y;k)=&D_{\tau}(x,y;k-1)\\
&+\sum_{d\ge1}x^dD_{\tau}^d(x,y;k-1)(D_{\tau}(x,y;k-1)-1)\\
&+\sum_{d\ge1}x^dD_{\tau}^{d-1}(x,y;k-1))\\
&+\sum_{d\ge1}x^dyD_{\tau}^{d-1}(x,y;k-1)(D_{\tau}(x,y;k-1)-1).
\end{split}
\]
Hence,
\[
\begin{split}
F_{\tau}(x,y;k)&=\frac{F_{\tau}(x,y;k-1)}{1-xD_{\tau}(x,y;k-1)},\\
D_{\tau}(x,y;k)&=\frac{(1-x+xy)D_{\tau}(x,y;k-1)+x(1-y)}{1-xD_{\tau}(x,y;k-1)}.
\end{split}
\]
Together with $F_{\tau}(x,y;0)=D_{\tau}(x,y;0)=1$,
$F_{\tau}(x,y;1)=D_{\tau}(x,y;1)=1/(1-x)$ and induction on $k$,
this yields the following result.
\begin{theorem}\label{the123}
Let $\tau=123$ be a subword pattern. For all $k\ge 2$, we have
\[
F_{\tau}(x,y;k)=\frac{1}{1-kx-\sum_{j=3}^k(-x)^j\binom{k}{j}(1-y)^{\lfloor
j/2\rfloor}U_{j-3}(y)},
\]
where $U_0(y)=U_1(y)=1$, $U_{2n}(y)=(1-y)U_{2n-1}(y)-U_{2n-2}(y)$,
and $U_{2n+1}(y)=U_{2n}(y)-U_{2n-1}(y)$. Furthermore, the
generating function for $U_n(y)$ is given by
\[
\sum_{n\ge 0}U_n(y)z^n=\frac{1+z+z^2}{1+(1+y)z^2+z^4}.
\]
\end{theorem}

Finally, Theorem \ref{thab} for $l=3$ and $m=2$ provides already the answer
for the last class. The corresponding result is summarized below.

\begin{theorem}\label{the213}
Let $\tau=213$ be a subword pattern. Then for all $k\ge 2$, we have
\[
F_{\tau}(x,y;k)=\frac{1}{1-x-x\sum_{i=0}^{k-2}\prod_{j=0}^i
(1-jx^2(1-y))}.
\]
\end{theorem}

\section{Further results}

We say that a subword pattern $\tau\in [m]^l$ is \emph{primitive}
if any two distinct occurrences of $\tau$ may overlap by at most
one letter. For example, the subword patterns $112$, $121$, $122$,
$132$, $211$, $212$, $213$, $221$, $231$, and $312$ are all the
primitive patterns of length three.

\begin{theorem}\label{thg2}
Let $\tau,\tau'\in[m]^l$ be two primitive subword patterns such
that there exists a permutation $\Phi\in S_l$ with $\Phi(1)=1$,
$\Phi(l)=l$ and $\tau'=\Phi\circ\tau$. In other words, $\tau$ and
$\tau'$ have the same supply of each letter, the same first letter
and the same last letter. Then $\tau$ and $\tau'$ are in the same
strong Wilf class.
\end{theorem}
\begin{proof}
Let $\sigma\in[k]^n$ contain $\tau$ exactly $r$ times. Since
$\tau$ is a primitive subword pattern, we can define a function
$f$ which changes any occurrence of $\tau$ in $\sigma$ to an
occurrence of $\tau'$. It is easy to see from the definition of 
primitive patterns that $f$ is a bijection, hence the theorem
follows.
\end{proof}

An immediate corollary is the following.

\begin{corollary}\label{the1322}
The subword patterns $1232$ and $1322$ are in the same strong Wilf
class.
\end{corollary}

\begin{theorem}\label{thg3}
All primitive subword patterns $\tau\in [m]^l$ such that
$\tau(1)=a$ and $\tau(l)=b$, where $a<b$, are in the same strong
Wilf class.
\end{theorem}
\begin{proof}
Similarly as in the proof of Theorem \ref{thll}, we get
\[
F_{\tau}(x,y;k)=F_{\tau}(x,y;k-1)+xF_{\tau}(x,y;k-1)D_{\tau}(x,y;k;a),
\]
where $D_\tau(x,y;k;h)$ is the generating function for the number
of words $\sigma\in [k]^n$ such that $(h,\sigma)$ contains $\tau$
exactly $r$ times.

Now let us consider the case $h=a+p(b-a)$. Let $\sigma=(\sigma',
h+b-a ,\sigma'')$, where $\sigma'$ is a word on the letters in
$[k]$ which does not contain $h+b-a$. Using the fact that $\tau$
is a primitive subword pattern, we get
\[
\begin{split}
&D_\tau(x,y;k;h)=D_\tau(x,y;k-1;h)\\
&+x^{l-1}y\binom{h-1}{a-1}\binom{k-(h+b-a)}{m-b}D_\tau(x,y;k;h+b-a)\\
&+x\left[F_\tau(x,y;k-1)-x^{l-2}\binom{h-1}{a-1}\binom{k-(h+b-a)}{m-b}\right]
D_\tau(x,y;k;h+b-a).
\end{split}
\]

Hence, by induction on $p$ and $k$, using
$F_\tau(x,y;m-1)=1/(1-(m-1)x)$ and $D_\tau(x,y;k;h)=0$ for $h>k$,
we get the desired result.
\end{proof}

Using the proof of the above theorem, we get the following
generalization.

\begin{corollary}\label{thg4}
Let $\tau,\tau'\in[m]^l$ be two primitive subword patterns such
that $\tau(1)=\tau'(1)=a$ and $\tau(l)=\tau'(l)=b$, where $a<b$.
Then the subword patterns
\[
\underbrace{aa\dots a}_p \tau\underbrace{b\dots bb}_p \quad \text{
and } \quad \underbrace{aa\dots a}_p \tau' \underbrace{b\dots
bb}_p
\]
are in the same strong Wilf class.
\end{corollary}

Theorem \ref{thg3} implies as well the following corollary.
\begin{corollary}\label{the1332} \
\begin{enumerate}
\item The subword patterns $1132$, $1232$, $1322$, and $1332$ are in the
same strong Wilf class.
\item The subword patterns $1432$ and $1342$ are in the same strong Wilf
class.
\end{enumerate}
\end{corollary}

\begin{theorem}\label{thg5}
Let $12\tau1\in[m]^{l+3}$ be a primitive subword pattern, and let
$\tau'$ be the same pattern $\tau$ with $1$ replaced by $2$. Then
the subword patterns $12\tau13$ and $12\tau'23$ are in the same
strong Wilf class.
\end{theorem}
\begin{proof}
If $\sigma\in[k]^n$ contains $12\tau13$ exactly $r$ times, then we
define $\sigma'$ as follows. If $(\sigma_i,\dots,\sigma_{i+l+3})$
is an occurrence of $12\tau13$, then we define
$\sigma'_{i+j}=\sigma_{i}+1$ for all $j$ such that
$\sigma_{i+j}=\sigma_i$. The function $f$ defined by
$f(\sigma)=\sigma'$ is a bijection since $12\tau1$ is a primitive
subword pattern.
\end{proof}

\begin{corollary}
The subword patterns $1213$ and $1223$ are in the same strong Wilf
class.
\end{corollary}

\small{
\begin{center}
\textsc{Acknowledgement}
\end{center}

The authors would like to thank the anonymous referee for bringing
several references, particularly \cite{RS}, to their attention.
The final version of this paper was written
while the second author (T.M.) was visiting University of Haifa, Israel in
January 2003. He thanks the HIACS Research Center and the Caesarea
Edmond Benjamin de Rothschild Foundation Institute for
Interdisciplinary Applications of Computer Science for financial
support, and professor Alek Vainshtein for his generosity. }


\end{document}